\documentclass[11pt,letterpaper,reqno]{amsart}
\allowdisplaybreaks 
\belowdisplayskip
\abovedisplayskip

\numberwithin{equation}{section} 

\usepackage{pgfplots}
\pgfplotsset{compat=1.18}  
\usetikzlibrary{intersections} 
\usepgfplotslibrary{fillbetween}

\usepackage{lmodern}
\usepackage[T1]{fontenc}

\makeatletter
\@namedef{subjclassname@2020}{%
  \textup{2020} Mathematics Subject Classification}
\makeatother
\usepackage[table,xcdraw]{xcolor}
\usepackage[a4paper,top=2.5cm,bottom=2.5 cm,left=2 cm,right=2 cm]{geometry}

\usepackage[utf8]{inputenc} 
\usepackage[english]{babel} 
\usepackage{changes}
\usepackage{fancyhdr}
\usepackage{amssymb}
\usepackage[toc,page]{appendix}
\usepackage{amsmath,amssymb,mathrsfs,mathrsfs,pifont,cancel,paralist,mathbbol,bbm}
\usepackage{amsthm}
\usepackage{mathtools}
\usepackage{caption, subcaption}
\usepackage{tabularx}
\usepackage{graphicx}
\usepackage{mathrsfs}
\usepackage{pdfcomment}
\usepackage{url}
\hypersetup{hidelinks}
\usepackage[a-1b]{pdfx}
\usepackage[T1]{fontenc}
\usepackage{bm}
\usepackage{enumitem}
\usepackage{comment}
\usepackage[all]{xy}
\usepackage{tikz}
\usepackage{tikz-cd}
\usepackage{mathbbol}
\usepackage{epigraph}
\usepackage{setspace}
\usepackage{ stmaryrd }
\usepackage{nicefrac,xfrac}
\usepackage{amsmath}
\usepackage{makecell} 
\usepackage{hyperref}

\definecolor{forestgreen(web)}{rgb}{0.13, 0.55, 0.13}

\newcommand{\commentoo}[1]{}
\newcommand{\N}{\mathbb{N}}
\newcommand{\Z}{\mathbb{Z}}

\newcommand{\F}{\mathbb{F}}

\newcommand{\Q}{\mathbb{Q}}

\newcommand{\R}{\mathbb{R}}
\newcommand{\C}{\mathbb{C}}
\renewcommand{\P}{\mathbb{P}}
\newcommand{\D}{\mathbb{D}}

\renewcommand{\d}{\mathrm{d}}
\newcommand{\mcal}{\mathcal}
\newcommand{\mbf}{\mathbf}
\newcommand{\summ}[1]{\sum_{\substack{#1}}}

\newcommand{\pol}[2]{\mathrm{Pol}_{#1}(#2)}

\newcommand{\ord}[1]{{\left| {#1}\right|}}

\newcommand{\e}{\mathrm{e}}

\newcommand{\sfree}{\text{ square-free}}
\newcommand{\leg}[2]{\left(\frac{#1}{#2}\right)}
\theoremstyle{definition}
\newtheorem{defi}{Definition}[section]

\newtheorem{remark}[defi]{Remark}

\newtheorem*{notation}{Notation}
\newtheorem*{outline}{Outline of the paper}
\newtheorem*{ai}{AI disclosure statement}

\newtheorem*{ack}{Acknowledgements}

\theoremstyle{plain}
\newtheorem{prop}[defi]{Proposition}
\newtheorem{coro}[defi]{Corollary}
\newtheorem{lemma}[defi]{Lemma}
\newtheorem{theorem}[defi]{Theorem}

\title{Ranks of fibers in random elliptic surfaces}
\author{Giacomo Bortolussi}
\address{Università di Roma Tor Vergata\\ Dipartimento di Matematica\\ 00133\\ Rome}
\email{bortolussi@axp.mat.uniroma2.it}
\begin{document}
\begin{abstract}
    For random integer polynomials $P$ of arbitrary degree we prove that for $100\%$ of elliptic surfaces $y^2=x^3-P(t)^2x$ the average 2-Selmer rank of the fibers is 3. 
    Consequently, for $100\%$ these random elliptic surfaces we establish a lower bound for the proportion of fibers attaining algebraic rank $0$ and $1$.
\end{abstract}
\maketitle

\setcounter{tocdepth}{1}

\section{Introduction}
For any polynomial \(P\in \mathbb Z[t]\), one may study the ranks over $\mathbb Q$ of the fibers of the elliptic surface
\begin{equation}\label{es}
E_{P(t)}\colon \qquad y^2=x^3-P(t)^2x.
\end{equation}
This is motivated by Silverman's conjecture, which concerns the statistical behavior of ranks of fibers of general elliptic surfaces.  Silverman's conjecture remains open in broad generality.
Given an elliptic curve $E/\mathbb Q$, it is in general difficult to obtain information about its algebraic rank $r(E)$. A more tractable invariant is its $2$-Selmer rank 
$s(E)\coloneqq \dim_{\mathbb F_2}\mathrm{Sel}_2(E/\Q)-\dim_{\mathbb F_2}E(\Q)[2]$,
which satisfies the inequality $r(E)\le s(E)$.
For a fixed $P\in \Z[t]$, it is natural to study for the \textit{polynomial twist family} (\ref{es}) the average 
$$\frac{1}{X}\sum_{\substack{1\le m\le X\\ P(m)\neq 0}}
2^{s(E_{P(m)})},$$
In general, such a quantity is expected to depend on the geometry of the specific (isotrivial) elliptic surface.

For $H\ge1$ and $d\ge 1$, define $$\pol{d}{H}\coloneqq\{a_dt^d+\ldots+a_1t+a_0\in \Z[t]: a_0,\ldots,a_d\in (0,H]\}\,.$$
In this work, we study, for \textit{random elliptic surfaces} of the form \eqref{es}, the average normalized $2$-Selmer cardinality of the fibers $2^{s(E_{P(m)})}$ in the \textit{exponential range} (\ref{exp_range}).
\begin{theorem}\label{teo_mainES}
    Fix $d\ge 1$, $\delta>\delta_0=48+(d+1)/2$, $\epsilon>0$ and $b\in \{1,\ldots,8\}$. Let $H>3$ and let $X>1$ satisfy \begin{equation}\label{exp_range}(\log H)^{\delta_0}\le X\le (\log H)^{\delta}\,.\end{equation} Then we have $$\frac{1}{\#\pol{d}{H}}\sum_{P\in \pol{d}{H}}\Bigg|\sum_{\substack{m\le X\\P(m)\equiv b\bmod 8}}\left(2^{s(E_{P(m)})}-3\right)\Bigg|^2\ll \frac{X^2}{(\log H)^{7/16-\epsilon}}\,,$$ where the implied constant depends only on $d,\delta,\epsilon$.
\end{theorem}
 By directly applying Chebyshev's inequality, we have the following result.
\begin{coro}\label{coro100}
    Keep the setting of Theorem \ref{teo_mainES}. For $100\%$ of polynomials $P\in \pol d H$ we have
    $$ \sum_{\substack{m\le X\\P(m)\equiv b\bmod 8}}2^{s(E_{P(m)})}= 3\#\{m\le X: P(m)\equiv b\bmod 8\}+o\bigl(X\bigr)\,.$$
   Fix $\epsilon>0$. For $100\%$ of $P\in \pol d H$ we have that
    $$ {\#\left\{m\le X\mid P(m)\equiv b\bmod 8\textrm{ and }r(E_{P(m)})\le 1\right\}}\ge\left(\frac 1 3 -\epsilon\right){\#\left\{m\le X(H)\mid P(m)\equiv b\bmod 8\right\}}\,.$$
   Assuming the parity conjecture, then for $100\%$ of $P\in \pol d H$ we have
    \begin{align*}&{\#\left\{m\le X\mid P(m)\equiv 1,2,3\bmod 8,\,r\bigl(E_{P(m)}\bigr)=0\right\}}\ge\left(\frac 1 3-\epsilon\right) {\#\left\{m\le X\mid P(m)\equiv 1,2,3\bmod 8\right\}}\end{align*}
    and 
    \begin{align*}
    {\#\left\{m\le X\mid P(m)\equiv 5,6,7\bmod 8,\,r(E_{P(m)})=1\right\}}\ge\left(\frac 5 6-\epsilon\right) {\#\left\{m\le X\mid P(m)\equiv 5,6,7\bmod 8\right\}}\,.\end{align*}
\end{coro}

Our polynomial-twist family  $E_{P(t)}$(\ref{es}) is a thin subfamily of the congruent number family  \begin{equation}
    E_D:\qquad y^2=x^3-D^2x\,.\label{cnf}
\end{equation}Goldfeld conjecture for $E_D$ predicts rank $0$ for $100\%$ of $D\equiv 1,2,3\bmod 8$ and  rank $1$ for  $100\%$ of $D\equiv 5,6,7\bmod 8$.  
For a fixed $P\in \Z[t]$ with $\deg P\ge2$, the values attained by $P$ over $\N$ have zero density in $\N$, thus Goldfeld conjecture gives no prediction for $E_{P(t)}$.
Corollary \ref{coro100} provides evidence supporting the minimalist philosophy for the \textit{typical} elliptic surface of the form (\ref{es}).

Although we have a fairly detailed conjectural understanding of the ranks in families of elliptic curves, unconditional information on the ranks of individual fibers remains difficult to obtain.
For certain isotrivial families, the situation is more tractable through the study of Selmer groups.
Recent breakthroughs by Smith use the structure of $2^\infty$-Selmer groups to prove that $100\%$ of quadratic twists of a given elliptic curve over $\Q$ have rank $\le 1$; combined with the Birch--Swinnerton--Dyer conjecture, this yields Goldfeld's conjecture \cite{Sm}. 

The study of $2$-Selmer statistics was initiated by Heath-Brown \cite{HB,HB2}, who built a $2$-descent counting machinery to study the average $2$-Selmer cardinality in the congruent number family. He proved that for  $b\in \{1,3,5,7\}$ one has \begin{equation}\label{teo_HBorig}\sum_{\substack{D\le T\\D\equiv b\bmod 8\\D\sfree}}2^{s(E_D)}=3\#\left\{D\le T\mid D\equiv b\bmod 8,\, D\sfree\right\}+O\left(\frac{T(\log\log T)^8}{(\log T)^{1/4}}\right)\,.\end{equation}
The same result is extended to $b\in \{2,6\}$ in \cite{Bo2}. In fact, in \cite{Bo2} we proved 2-Selmer cardinality asymptotics for (\ref{cnf}) on arbitrary arithmetic progressions $D\equiv b\bmod q$ with $q$ allowed to go to infinity slowly with $T$.
\begin{remark}
    Although in the congruent number family (\ref{cnf}) we expect the algebraic rank of the fibers to depend on the congruence class of $D$ modulo $8$, the same bias is not present for the average $2$-Selmer cardinality. To spot modulo $8$ biases, one needs to look at the distribution of the $2$-Selmer ranks \cite{HB2}.
    \end{remark}

 Heath-Brown's work was subsequently extended in several directions; see, for example, Swinnerton-Dyer \cite{SD}, Kane \cite{Ka} and Klagsbrun--Mazur--Rubin \cite{KMR}.
By completely different methods, Bhargava--Shankar \cite{BS} proved that the average of $2^{s(E_{a,b})}$ in the family $\{E_{a,b}\}_{a,b\in \Z}\coloneqq\{y^2=x^3-ax-b\}_{a,b\in \Z}$ of all elliptic curves is $3$.

Already Koymans--Pagano--Sofos \cite{KPS} studied polynomial twist families of elliptic curves with full rational $2$-torsion, and obtained matching  upper and lower bounds for the order of magnitude of the general moments of the $2$-Selmer cardinality. 
Later, Koymans--Smith \cite{KS} conjectured the correct order of magnitude for 2-Selmer moments in generic elliptic families.
Motivated by \cite[Conjecture 1]{KS}, we provide evidence for the leading constant of the first moment for the typical elliptic surface of the form (\ref{es}).

\begin{outline}
    In \S \ref{sec_prel} we give some preliminaries on the Selmer rank and we recall the necessary notions from \cite{Bo2}. Relying on Frei--Sofos' \cite{FS} heat kernel machinery, in \S \ref{section-secon_moment} we show that an upper bound for the second moment over random polynomials follows from uniform equidistribution on small moduli arithmetic progressions, which is provided in \cite{Bo2} for type 2--4 error terms. In \S \ref{sec-g3} we separately take care of the second moment of type 1 and R error terms using Nair's bound \cite{Na}.  In \S \ref{sect-proof} we state our main result on square-free polynomial twists, namely Theorem \ref{teo_main}.
    Most of the arguments are carried out with square-free twisting parameters: in \S \ref{sec-ssf} we use the square-free sieve to remove the square-free condition, i.e. we show that  Theorem \ref{teo_mainES} follows from Theorem \ref{teo_main}.
\end{outline}

\begin{notation}
    From now on, $D$ will always be an odd, square-free positive integer. We use $\tau_k(n)$ for the number of decompositions of $n$ as a product of exactly $k$ positive integers; $\tau_2=\tau$ is the classical divisor-counting function. We use $\omega(n)$ for the number of distinct prime divisors of $n$. Both the algebraic rank and the $2$-Selmer rank are always intended over $\Q$.
\end{notation}

\begin{ai}
     Every idea in the present work is human. AI was used only to improve exposition and style.
\end{ai}
\begin{ack}
    The author would like to thank his supervisor, Efthymios Sofos, for introducing him to the topic and for the support during the writing process. 
\end{ack}

\tableofcontents
\section{Preliminaries}\label{sec_prel}
Let $E$ be an elliptic curve over $\Q$. By Mordell's theorem, its group of rational points is finitely generated and admits the decomposition
$$E(\Q)\cong E_{\mathrm{tors}}\oplus \Z^{r(E)}\,.$$
Although the torsion subgroup is well understood, determining the algebraic rank $r(E)$ is typically a difficult problem. A standard approach to bounding the rank relies on Selmer groups. Indeed, the natural injection
$$\frac{E(\Q)}{2E(\Q)}\hookrightarrow\mathrm{Sel}_2(E/\Q)$$
implies that $2^{r(E)}\cdot\#E(\Q)[2]\leq \#\mathrm{Sel}_2(E/\Q)$. Defining $s(E)\coloneqq \dim_{\F_2}\mathrm{Sel}_2(E/\Q)-\dim_{\F_2}E(\Q)[2]$, we obtain the upper bound $r(E)\leq s(E)$. Abusing notation, we write $s(D)$ and $r(D)$ for $s(E_D)$ and $r(E_D)$. Note that our family $E_P$ has full 2-torsion, so $\#E_{P(m)}(\Q)[2]=4$.

By explicitly counting $2$-descents, Heath-Brown obtained an explicit counting formula for $2^{s(D)}$ when $D$ is odd and square-free \cite[Lemma 3]{HB}. Relying on Monsky's matrix approach, the author proved an analogous formula for $2^{s(2D)}$ \cite[Formula 2.7]{Bo2}, again with $D$ odd and square-free.
As an immediate consequence of the two formulae, one has \begin{equation}\label{bound_sel}
    2^{s(D)}\le 4^{\omega(D)}\qquad\text{and}\qquad 2^{s(2D)}\le 4^{\omega(D)}\,.
\end{equation}

Fix three integer parameters $z,C,h>0$. Then, depending on this choice, for any positive odd square-free integer $D$ we have \begin{equation}
    \label{eq-2selcard}
    2^{s(D)}=3+\Bigg[G^{(1)}+\summ{\mbf A \text{ multi-interval}\\\text{of type 2}}G^{\mbf A}+\summ{\mbf A \text{ multi-interval}\\\text{of type 3}}G_h^{\mbf A}+\sum_{\mcal J\in \mathscr B}\left(G^{(\mcal J,\text{a})}+G_h^{(\mcal J,\text{b})}\right)-G^{(\text{R})}\Bigg](D)\,.
\end{equation}
This formula displays Heath-Brown's error-term decomposition, refined as in \cite[Equation (3.8)]{Bo2}.
In the formula above, the constant term $3$ contributes to the main term, while each $G^{\Omega}_h(D)$, which may depend on $C,z,h,\mbf A$, and contributes to an error term. The two sums over $\mbf A$ run through $\ll (\log D)^{16}$ dyadic multi-intervals, while the sum over $\mcal J$ runs through the 24 possible four-tuples of unlinked variables as identified in \cite[Lemma 9]{HB} or \cite[Lemma 2.1]{Bo2}. In particular, restricting to $D\le T $ one needs to consider $\ll (\log T)^{16}$ dyadic multi-intervals. To keep the exposition as self-contained as possible, we avoid further details on the different types of error terms. The reader should keep in mind that:
\begin{itemize}
    \item For type 2--4 error terms (i.e. $G^\mbf A$ and $G^{(\mcal J,-)}$), \cite[Propositions 3.6, 3.8, 3.14, 3.16]{Bo2} proves uniform equidistribution in arithmetic progressions of modulus $q\le (\log T)^N$. In particular, the saving is strong enough to allow an application of Frei--Sofos' heat-kernel machinery.
    \item For type 1 and R error terms, the saving in arithmetic progressions is too weak to allow an application of Frei--Sofos' heat-kernel machinery; in these two cases we instead open up the second moment explicitly and employ Nair's bound \cite{Na} to obtain a satisfactory upper bound.
\end{itemize}
Furthermore, for each $h$ and each $G_h^\Omega$ appearing in (\ref{eq-2selcard}) we have \begin{equation}
    \label{eq_G_Omega}\ord{G_h^\Omega(D)}\ll 4^{\omega(D)}\,.
\end{equation}
A formula similar to (\ref{eq-2selcard}) is available for $2^{s(2D)}$, and the corresponding versions of the relevant propositions in \cite{Bo2} also hold true.
\begin{remark}
    In \cite{Bo2} the established error term is uniform, for every $N>0$, in $q\le (\log T)^N$; the same result is no longer true for arbitrary arithmetic progressions with $q\le T^\epsilon$ for some $\epsilon>0$, see \cite[Example 3.19]{Bo2}. Since uniform saving in  arithmetic progressions is available only for $q\le (\log T)^N$, in the current work we need to adapt Frei--Sofos' machinery so as to fit into the exponential range $X(H)\sim(\log H)^{\delta}$.
\end{remark}
\begin{remark}
    Maintaining the dependence on $h$ in each $G^{\Omega}_h$ allows one to prove uniform equidistribution in arithmetic progressions. A remarkable feature of \cite{Bo2} is that, even though the functions $G_h^\Omega$ may depend on $h$, the implied constants in the quoted bounds are uniform in $h$.
\end{remark}
 \begin{remark}
     As $E_D\cong E_{k^2D}$ for any $k\in \Z_{>0}$, it is natural to work with square-free twisting parameters, as in  (\ref{teo_HBorig}). However, it would feel unnatural to remove square-full fibers in Theorem \ref{teo_mainES}, mostly because the density of square-free fibers depends on the specific $P$. For most of the paper we work with square-free values (see Theorem \ref{teo_main}), as this is the correct setting for Heath-Brown's $2$-descent counting formula. In particular, each $G^\Omega$ appearing in (\ref{eq-2selcard}) is by definition supported on square-free odd integers (for both $2^{s(D)}$ and $2^{s(2D)}$). In \S\ref{sec-ssf} we use the square-free sieve to remove the square-free restriction, thereby deducing Theorem \ref{teo_mainES} from Theorem \ref{teo_main}.
 \end{remark}

\section{Main theorem and strategy of the proof}\label{sect-proof}
For $P\in \pol d H$, $k\in \Z_{>0}$, and $b\in \{1,\ldots,8\}$, define
\begin{equation}\label{defpsi}
\Psi_{P,k,b}(X)\coloneqq \sum_{\substack{m\le X\\k^2\mid P(m)\\P(m)/k^2\sfree\\P(m)/k^2\equiv b\bmod 8}}\left(2^{s\left(\frac{P(m)}{k^2}\right)}-3\right)\,.
\end{equation}

The main result of our work is the following.
\begin{theorem}
    \label{teo_main}
    Let $d\ge 1, \nu,\epsilon>0$, $\delta>\delta_0=48+{(d+1)/2}$ and $b\in \{1,2,3,5,6,7\}$. Let $H>3$ and let $X=X(H)$ satisfy $(\log H)^{\delta_0}\le X\le (\log H)^{\delta}$. Then for any $1\le k\le (\log H)^{\nu}$ we have
    \begin{equation}\label{eq-main}
    \frac{1}{\#\pol d H}\sum_{P\in \pol d H}\Big|\Psi_{P,k,b}(X)\Big|^2\ll \frac{X^2}{k^2(\log H)^{7/16-\epsilon}}\,.
    \end{equation}
    where the implied constant depends only on $d,\delta,\nu,\epsilon$.
\end{theorem}

\begin{remark}
    The divisibility condition imposed in Theorem \ref{teo_main} is new in the context of random polynomials, and it is crucial for us since it allows us to sum also over square-full polynomial twists (see \S\ref{sec-ssf}). Special care is taken in \S\ref{section-secon_moment} to record the dependence of the error terms on $k$.
\end{remark}

\begin{remark}
    In general, $\#\pol d H\sim H^{d+1}$. The density of square-free values assumed by $P$ is influenced by its local behavior. Although in general very little is known for fixed polynomials, Sofos \cite{So} proved that for $100\%$ of $P\in \pol d H$ the density is the expected one for $x$ in the range $x(H)\ll H^{\alpha}$.
\end{remark}

Before starting the actual proof, we need a technical tool to separate the contributions of the different types of terms in the second moment.

\begin{lemma}
    \label{lemma-cs}  Let $f_1,\ldots,f_n\colon \N\to \R$. Then for every $H,X>1$ we have $$ {\frac{1}{\#\pol d H}\!\sum_{P\in \pol d H}\!\!\left(\sum_{m\le X}f_1(P(m))+\ldots+f_n(P(m))\!\right)^2}\!\!\!\le {\frac n {\#\pol d H}\sum_{i=1}^n\sum_{P\in \pol d H}\!\!\left(\sum_{m\le X}f_i(P(m))\!\!\right)^2}\!.$$
\end{lemma}
\begin{proof}
    Apply the Cauchy--Schwarz inequality and the triangle inequality:
    \begin{align*}&  {\mathbb E_{\pol d H}\Bigg[\Bigg(\sum_{m\le X}\sum_{1\le i\le n}f_i(P(m))\Bigg)^2\Bigg]}
    = {\mathbb E_{\pol d H}\Bigg[\sum_{m_1,m_2\le X}\sum_{1\le i,j\le n}f_i(P(m_1))f_j(P(m_2))\Bigg]} \\
    &\le
    \sum_{1\le i,j\le n} {\bigg|\mathbb E_{\pol d H}\Bigg[\Bigg(\sum_{m_1\le X}f_i(P(m_1))\Bigg)\Bigg(\sum_{m_2\le X}f_j(P(m_2))\Bigg)\Bigg]\Bigg|}
    \\&\le
    \sum_{1\le i,j\le n} {\mathbb E_{\pol d H}\Bigg[\Bigg(\sum_{m_1\le X}f_i(P(m_1))\Bigg)^2\Bigg]}^{1/2} {\mathbb E_{\pol d H}\Bigg[\Bigg(\sum_{m_2\le X}f_j(P(m_2))\Bigg)^2\Bigg]}^{1/2}
    \\&=
    \Bigg(\sum_{i=1}^n\Bigg[ \mathbb E_{\pol d H}\Bigg(\sum_{m\le X}f_i(P(m))\Bigg)^2\Bigg]^{1/2}\Bigg)^2\,. \end{align*}
    Conclude by applying the Cauchy--Schwarz inequality in the form $\left(\sum_{i=1}^n\sqrt{a_i}\right)^2\le n\sum _{i=1}^n a_i$.
\end{proof}

\begin{proof}[Proof of Theorem \ref{teo_main}]
We give a proof for $b\in \{1,3,5,7\}$, which allows one to use formula (\ref{eq-2selcard}) directly; the same strategy applies in the cases $b\in \{2,6\}$, but with $g_2$ in place of $g$ in (\ref{eq-2selcard}) (see \cite[\S\S 2.4, 3.5]{Bo2}).

For some $B,\eta>0$ to be chosen later, set
\begin{equation}\label{def_zC}
z(H)\coloneqq 2^{\big\lfloor \log_2\big((\log H)^B\big)\big\rfloor}
\qquad \text{and}\qquad
C(H)\coloneqq 2^{\big\lfloor (\log H)^{\eta}/\log 2\big\rfloor}\,.
\end{equation}
Then $(\log H)^B/2\le z\le (\log H)^B$ and $e^{(\log H)^{\eta}}/2\le C\le e^{(\log H)^{\eta}}$.

Using Lemma \ref{lemma-cs}, we obtain
\begin{align*}
    &\summ{P\in \pol d H}\Bigg|\summ{m\le X\\k^2\mid P(m)\\P(m)/k^2\equiv b\bmod 8\\P(m)/k^2\sfree}\left(2^{s\left(\frac{P(m)}{k^2}\right)}-3\right)\Bigg|^2
    \\&\ll \summ{P\in \pol d H}\Bigg|\summ{m\le X}^*G^{(1)}\bigg(\frac{P(m)}{k^2}\bigg)\Bigg|^2
    +\summ{P\in \pol d H}\Bigg|\summ{m\le X}^*G^{(\mathrm{R})}\bigg(\frac{P(m)}{k^2}\bigg)\Bigg|^2
    \\&\quad
    +\summ{P\in \pol d H}\Bigg|\summ{m\le X}^{*}\Bigg(2^{s\left(\frac{P(m)}{k^2}\right)}-3-G^{(1)}\bigg(\frac{P(m)}{k^2}\bigg)+G^{(\mathrm{R})}\bigg(\frac{P(m)}{k^2}\bigg)\Bigg)\Bigg|^2\,.
\end{align*}
where, for brevity, $(*)$ indicates a repetition of the conditions on $m\le X$ made explicit in the first sum.
The two second moments involving $G^{(1)}$ and $G^{(\mathrm{R})}$ are bounded in Propositions \ref{prop_sm3} and \ref{prop_smR}, respectively, and the two resulting error terms represent the main error terms in Theorem \ref{teo_main} with $0<\eta<\min\left\{1/{16},{8\epsilon}/{377}\right\}$.

As for the remaining part, from (\ref{eq-2selcard}) we define
\begin{equation*}
    F\coloneqq \mu^2(\cdot)\left(2^{s(\cdot)}-3\right)-G^{(1)}+G^{(\mathrm{R})}=\summ{\Omega}G_h^{\Omega}\,,
\end{equation*}
where $h$ is still to be specified and $\Omega$ runs through the $\ll (\log H)^{16}$ indices remaining from (\ref{eq-2selcard}) after removing $(1)$ and $(\mathrm{R}) $.
Note that $h$ does not appear in every $G_h^{\Omega}$; we nevertheless keep this (sometimes unnecessary) index to avoid splitting into cases. Recall also that the functions $G_h^{\Omega}$ are, by definition, supported on square-free integers.

Calling $F_b\coloneqq F\cdot \mathbb 1_{\equiv b\bmod 8}$, we are left to bound the second moment of $F_b$ over random polynomials. By repeating the argument in the proof of Theorem \ref{theorem-second_moment}, 
$$\summ{P\in \pol d H}\Bigg|\summ{m\le X}^*F\bigg(\frac{P(m)}{k^2}\bigg)\Bigg|^2\ll \summ{m\le X}\ord{\mbf J_{F_b, F_b,k}(m,m)}+\summ{m_1\neq m_2\le X }\ord{\mbf J_{F_b, F_b,k}(m_1,m_2)}\,.$$
    From (\ref{eq-2selcard}) and (\ref{bound_sel}), we deduce $F_b(D)\ll 4^{\omega(D)}$, which, combined with Proposition \ref{prop-diag2}, gives for the diagonal part above $$\ll \frac{H^{d+1}X(\log H)^{15+(d+1)/2}}{k^2}\,.$$

    Concerning the non-diagonal part, fix $m_1,m_2$, let  $q(m_1,m_2)\coloneqq \mathrm{lcm}(\ord{m_1-m_2},8)$  and let $h=h(m_1,m_2)$ be the square-full part of $q$ (as in \cite[\S 3]{Bo2}). Note that $q,h\ll (\log H)^\delta$. 
    Calling $\tilde G^\Omega_h\coloneqq G^\Omega_h\cdot \mathbb 1_{\equiv b\bmod 8}\cdot \mathbb 1_{\le (d+1)HX^d/k^2}$, then $$J_{F_b, F_b,k}, (m_1,m_2)=\summ{\Omega_1,\Omega_2}\mbf J_{\tilde G^{\Omega_1}_h,  {\tilde G^{\Omega_2}_h},k}(m_1,m_2)$$
    The choice of $q$ is made to include the further condition modulo 8, in fact to apply Proposition \ref{teo_fsrevisited}, we need uniform savings for sums of the form
\begin{equation}\label{eq_G_ap}
    \summ{D\le T\\D\equiv a\bmod q}G_h^{\Omega_i}(D)
\end{equation}
over all arithmetic progressions $\{qn+a\}$ with $a\equiv b\bmod 8$.
For each $\Omega$ appearing in $F$, this is achieved in \cite{Bo2} according to Table \ref{table_props}, with the following straightforward modification: define $z$ and $C$ to depend on $H$ (and not on $T$) according to (\ref{def_zC}); then the bounds in Table \ref{table_props} hold as long as $T\gg H^\epsilon$ for some small $\epsilon$.
\begin{equation}\label{table_props}
    \begin{tabular}{|c|l|c|}\hline
    Type of $\Omega$     & (\ref{eq_G_ap}) bounded by & given by  \\
    \hline
    2 &$\ll \frac{q^2T(\log T)^{9/2}}{z^{1/4}}$ & \cite[Proposition 3.6]{Bo2}\\\hline
    3 & $\ll \frac{T(\log T)^{8}}{e^{c(\log H)^{\eta/2}}}$& \cite[Proposition 3.8]{Bo2} \\\hline
    4a/4b & $\ll \frac{T(\log T)^2}{e^{c(\log H)^{\eta/2}}}$& \cite[Propositions 3.14, 3.16]{Bo2}\\\hline
    \end{tabular}
\end{equation}
We remark that, even though the different $G_h^{\Omega}$ may depend on $h$, the saving one obtains in every residue class is independent of $h$.

Since for every choice of $B,\eta>0$ the asymptotically larger term always comes from $\Omega$ of type 2, we may assume the quantity (\ref{eq_G_ap}) to be bounded by $\ll \frac{q^2T(\log T)^{9/2}}{z^{1/4}}$. Then Proposition \ref{teo_fsrevisited} applies with $f_i=G_h^{\Omega_i}\cdot \mathbb 1_{\equiv b\bmod 8}$, and for fixed $m_1,m_2,\Omega_1,\Omega_2$ we have
$$\ord{\mbf J_{\tilde G^{\Omega_1}_h,  {\tilde G^{\Omega_2}_h},k}(m_1,m_2)}\ll \frac{H^{d+1}X^{4d+6}(\log X)^2(\log H)^9}{k^4z^{1/2}}\,.$$
There are $\ll (\log H)^{16}$ choices for each $\Omega_i$, so the non-diagonal part in the second moment of $F$ is bounded by
$$\ll \frac{H^{d+1}X^{4d+8}(\log H)^{41}(\log X)^2}{k^4z^{1/2}}\,.$$
Since $X\le (\log H)^\delta$, it is enough to take $B>1000d\delta$ to have a small enough error term.
\end{proof}


\section{Second moment and equidistribution}\label{section-secon_moment}
In this section, we prove that uniform bounds in small-modulus arithmetic progressions, together with a reasonable bound for a function $f$, allow one to bound the second moment (over random polynomials) of $\sum_{m\le X} f(P(m))$, with $X$ in the exponential range $\e^{(\log\log H)^{\frac {51}{100}}}\le X\le (\log H)^\delta$. The diagonal part is controlled using Nair's bound \cite{Na}, while the non-diagonal part is handled via an application of the circle method with heat kernels.

In the setting $X\le H^{c}$, Frei--Sofos \cite[Theorem 1.6]{FS} introduced heat kernels to bound the non-diagonal term of the second moment in the space of random binary forms. Using the same method, Chatzakos--Sofos \cite[Theorem 4.1]{CS} produced a similar result in the space of random polynomials, again with $X\le H^{c}$. Both results rely on some form of uniform equidistribution in large-modulus arithmetic progressions. 
In our case, uniform equidistribution is available only for small-modulus arithmetic progressions, placing us in the setting of \cite[Remark 1.7]{FS}.
 Here we adapt Chatzakos--Sofos' proof to the exponential range with the additional divisibility condition $k^2\mid P(m)$ for a fixed $k\le (\log H)^\nu$.
 
For $H\ge1$ and $\alpha\in\R$, define the heat kernel as $$K_H(\alpha)\coloneqq \sum _{n\in \Z}\exp\left(-\frac{\pi n^2}{H^2}+in\alpha\right)\,.$$
\begin{lemma}\cite[Lemma 2.15]{FS} \label{lemma-HK} Heat kernels have the following properties:\begin{enumerate}
    \item for every $\alpha\in \R$, $0\le K_H(\alpha)\ll H$;
    \item for every $c\in \Z$, $\widehat K_H(c)=e^{-\pi c^2/H^2}$ (where $\widehat K_H$ is the Fourier transform of $K_H$ in $[-\pi,\pi]$);
    \item for every $c\in \Z$, $\mathbb 1_{[0,H]}(c)\ll\hat K_H(c) $;
    \item $\int_{-\pi}^{+\pi} K_H(\alpha)\mathrm d\alpha=2\pi $;
    \item for every $\rho\in (0,\pi)$, $\int_{-\pi}^{+\pi}\mathbb 1_{|\alpha|>\rho} K_H(\alpha) \,\d \alpha\ll \frac{1}{\rho H\exp\left((\rho H)^2/(4\pi)\right)}$ as $H\to \infty$;
\end{enumerate}
\end{lemma}

\begin{theorem}
    \label{theorem-second_moment} Let $H>3,d\ge1$, let $\delta,\nu>0$, let $0<\epsilon <1/3$ and let $X=X(H)$ satisfy $\e^{(\log\log H)^{\frac {51}{100}}}\le X\le (\log H)^\delta$. Let $f\colon \Z\rightarrow \C$ be a function supported on positive integers such that, for some $u>0$, we have \begin{equation}\label{condfs1}\ord{f(n)}\ll u^{\omega(n)}\,.\end{equation}
    Suppose that there exists $A>0$ such that for all integers $q,a$ with $1\le q\le(\log H)^\delta$ we have \begin{equation}\label{condfs2}
       \Bigg| \sum_{\substack{n\le T\\n\equiv a \bmod  q}}f(n)\Bigg|\ll_A \frac {qT}{(\log T)^A}\qquad \qquad\text {for }T\ge H^\epsilon\,.
    \end{equation} Then, for every $1\le k\le (\log H)^\nu$, we have
    \begin{equation*}
        \frac 1 {\#\pol d H}\sum_{P\in \pol d H}\Bigg|\sum_{\substack{m\le X\\k^2\mid P(m)}}f\left(\frac{P(m)}{k^2}\right)\Bigg|^2\ll\frac{ X(\log H)^{u^2-1+(d+1)/2}}{k^2}+\frac{(\log X)^2X^{4d+6}}{k^4(\log H)^{2A}}\,,
    \end{equation*} where the implied constant depends only on $d,\delta,u,A,\nu,\epsilon$ and on the constant involved in (\ref{condfs1}) and (\ref{condfs2}).
\end{theorem}
\begin{proof}
Identify the set of integer polynomials of degree at most $d$ with $\Z^{d+1}$ by writing $P_{\mbf c}\in \Z[t]$ as $c_0+c_1t+\ldots +c_dt^d$. Replicating the idea of \cite[Equation 4.3]{CS}, we have
\begin{align}\notag
    \summ{P_{\mbf c}\in \pol d H}\Bigg|&\summ{m\le X\\k^2\mid P_\mbf c(m)}f\left(\frac{P_\mbf c(m)}{k^2}\right)\Bigg|^2\ll \sum_{\mbf c\in \Z^{d+1}}\Bigg(\prod_{j=0}^d\hat K_H(c_j)\Bigg)\Bigg|\summ{m\le X\\k^2\mid P_{\mbf c}(m)\\P_\mbf c(m)\le (d+1)X^dH}f\left(\frac{P_\mbf c(m)}{k^2}\right)\Bigg|^2\\
    &=\notag
    \summ{m_1,m_2\le X}\summ{\mbf c\in \Z^{d+1}\\k^2\mid P_\mbf c(m_1),\, P_\mbf c(m_2)\\P_\mbf c(m_1),P_\mbf c(m_2)\le (d+1)X^dH }\Bigg(\prod_{j=0}^d\hat K_H(c_j)\Bigg)  f\left(\frac{P_\mbf c(m_1)}{k^2}\right)\overline{f}\left(\frac{P_\mbf c(m_2)}{k^2}\right)\\
   &=
   \summ{m_1\neq m_2\le X} \mbf J_{f, f,k}(m_1,m_2)+\summ{m\le X}\mbf J_{f, f,k}(m,m)\,,\label{eq_interm}
\end{align}
where $\mbf J$ is defined as 
\begin{align}\label{def_J}
    \qquad\qquad \mbf J_{g_1,g_2,k}(m_1,m_2)\coloneqq & \frac 1 {(2\pi)^2} \int_{(-\pi,\pi]^2}{S_{g_1,k}(\alpha_1)}\overline{S_{g_2,k}(\alpha_2)}\left(\prod_{j=0}^dK_H(m_1^j\alpha_1-m_2^j\alpha_2)\right)\d\alpha_1\d\alpha_2
    \\\label{defS}&\text{ with } \qquad S_{g,k}(\gamma)\coloneqq \sum_{\substack{h\le (d+1)X^dH\\k^2\mid h}}g(h/k^2)e^{i\gamma h}.
\end{align}
The diagonal part of (\ref{eq_interm}) is bounded by Proposition \ref{prop-diag2}, while the non-diagonal part follows, after fixing $m_1,m_2$, from Proposition \ref{teo_fsrevisited} $\sigma=\tau=1,$ $W_H(t)=(\log (2+t))^A$.
\end{proof}

\subsection{Diagonal term} For the diagonal term, we exploit the freedom in the constant coefficient to impose the correct divisibility condition. Then, the bound (\ref{condfs1}) allows us to employ Nair's lemma \cite{Na}.

\begin{lemma}[Nair]\label{lemma_nair}
    Let $x,y>1$, $q,a\in \N^+$ with $\gcd(q,a)=1$ and $u>0$. Assume that for some $\alpha,\delta\in (0,1)$ and $c=c(\delta,\alpha)>0$ we have $x^{\alpha}\le y \le x$ and $q,a\le    cx^{\delta}$. Then $$ \sum_{x-y\le n\le x}u^{\omega(qn+a)}\ll C_{u,q}y (\log x)^{u-1}\,, $$ where the implied constant depends only on $u,\delta,\alpha$ and $C_{u,q}=\prod_{p\mid q}\frac{e^{-u/p}} {1-1/p}$.
\end{lemma}
\begin{proof}
    This follows directly from Nair's main theorem\cite{Na} with $f(n)=u^{\omega(n)}$ and $P(t)=qt+a$.
\end{proof}

Call $\pol d H'\coloneqq\{c_dt^d+\ldots+c_1t+c_0\in \Z[t]: \ord {c_i}\le H\text{ for }i=0,\ldots,d\}\supset \pol d H$.
\begin{prop}\label{prop-diag1}
Let $H>3,d\ge1$, let $\delta,\nu>0$ and let $X=X(H)$ satisfy $1< X\le (\log H)^\delta$. Let $f\colon \Z\rightarrow \C$ be a function supported on positive integers such that, for some $u>0$, we have \begin{equation}\label{cond_bound0}\ord{f(n)}\ll u^{\omega(n)}\,.\end{equation}
Then, for every $1\le k\le (\log H)^\nu$, we have \begin{equation*}
    \frac 1 {\#\pol d H'}\sum_{P\in \pol d H'}\sum_{\substack{m\le X\\k^2\mid P(m)}}\Bigg|f\left(\frac{P(m)}{k^2}\right)\Bigg|^2\ll \frac{X(\log H)^{u^2-1}}{k^2}\,,
\end{equation*}
    where the implied constant depends only on $d,\delta, u,\nu$ and on the constant involved in (\ref{cond_bound0}).
\end{prop}
\begin{proof}
    By assumption (\ref{cond_bound0}), the sum over $P$ is  \begin{align*}
        \ll\sum _{P\in \pol dH'}\sum_{\substack{m\le X\\k^2\mid P(m),\,P(m)>0}}u^{2\omega\left(\frac{P(m)}{k^2}\right)}\,.
    \end{align*}
    Writing $P(m)=c_0+c_1m+\ldots+c_dm^d\coloneqq c_0+N$, we isolate the coefficient $c_0$, set $c_0'\coloneqq\frac{N+c_0}{k^2}$, and consider only the integer values of $c_0'$; then the sum is
    $$\sum_{\substack{c_1,\ldots,c_d\in [-H,H]\\m\le X}}\sum_{\substack{c_0\in [-H,H]\\k^2\mid P(m),\,P(m)>0}}u^{2\omega\left(\frac{N+c_0}{k^2}\right)}=\sum_{\substack{c_1,\ldots,c_d\in [-H,H]\\m\le X}}\sum_{c_0'\in \left[\max\big\{0,\frac{N-H}{k^2}\big\},\frac{N+H}{k^2}\right]}u^{2\omega(c_0')}\,.$$
    By Lemma \ref{lemma_nair} (after enlarging, if necessary, to a positive interval of length $2H/k^2$) we can bound the sum above by \begin{equation*}
        \ll \sum_{\substack{c_1,\ldots,c_d\in [-H,H]\\m\le X}}\frac{H\left(\log H\right)^{u^2-1}}{k^2}\ll \frac{XH^{d+1}(\log H)^{u^2-1}}{k^2}\,.\qedhere\end{equation*}
\end{proof}

\begin{prop}\label{prop-diag2} Let $H>3,d\ge1$, let $\delta,\nu>0$ and let $X=X(H)$ satisfy $1< X\le (\log H)^\delta$. Let $f\colon \Z\rightarrow \C$ be a function supported on positive integers such that, for some $u>0$, we have \begin{equation}\label{cond_bound_2}\ord{f(n)}\ll u^{\omega(n)}\,.\end{equation}
Then, for every $1\le k\le (\log H)^\nu$, we have \begin{equation*}
    \frac{1}{\#\pol d H}\summ{m\le X}\mbf J_{f, f,k}(m,m)\ll \frac{X(\log H)^{u^2-1+(d+1)/2}}{k^2}\,,
\end{equation*}
    where the implied constant depends only on $d,\delta, u,\nu$\, and on the constant involved in (\ref{cond_bound_2}).
\end{prop}
\begin{proof}
     Let $\lambda=2d\sqrt{\log (XH)}$ and split $\sum_{m\le X}J_{f, f, k}(m,m)$ as \begin{equation}\label{eq_21}\summ{m\le X}\summ{\ord{\mbf c}>\lambda H\\k^2\mid P_\mbf c(m)\\P_\mbf c(m)\le (d+1)X^dH} \Bigg(\prod_{j=0}^d\hat K_H(c_j)\Bigg) \bigg|f\bigg(\frac{P_\mbf c(m)}{k^2}\bigg)\bigg|^2
    + \summ{m\le X}\summ{\ord{\mbf c}\le \lambda H\\k^2\mid P_\mbf c(m)\\P_\mbf c(m)\le (d+1)X^dH} \Bigg(\prod_{j=0}^d\hat K_H(c_j)\Bigg) \bigg|f\bigg(\frac{P_\mbf c(m)}{k^2}\bigg)\bigg|^2\,,\end{equation}
    where $\ord{\mbf c}=\max_j\ord {c_j}$. For any $m\le X$ bound  $\ord{f(P_\mbf c(m)/k^2)}^2\ll u^{2\omega(P_\mbf c(m)/k^2)}\ll \ord {\mbf c}X^d/k^2 $, then the first sum in (\ref{eq_21}) is bounded by \begin{align*}
        \ll \frac{X^{d+1}}{k^2}\summ{\ord{\mbf c}>\lambda H}\!\!\ord {\mbf c}\prod_{j=0}^de^{-\pi c_j^2/H^2}\ll \frac{X^{d+1}}{k^2}\left(\summ{c\in \Z}e^{-\pi c^2/H^2}\right)^d\!\!\left(\summ{c>\lambda H}e^{-\pi c^2/H^2}c\right)\ll \frac {X^{d+1}H^{d+3}}{k^2e^{\pi \lambda ^2}}\ll \frac{XH^{d+1}}{k^2}\,.
    \end{align*} 
    Use $\hat K_H(c_j)\ll 1$ together with Proposition \ref{prop-diag1} to bound the second sum in (\ref{eq_21}) by \begin{equation*}\ll \summ{P\in \pol d {\lambda H}'}\summ{m\le X\\k^2\mid P(m)}\bigg|f\bigg(\frac{P(m)}{k^2}\bigg)\bigg|^2\ll\frac{XH^{d+1}(\log H)^{u^2-1+(d+1)/2}}{k^2}\,.\qedhere\end{equation*}
    
\end{proof}

\subsection{Non-diagonal term}

\begin{prop}
    \label{teo_fsrevisited} Fix $d\ge 1$, $\delta,\nu>0$. Let $H>3$, let $X=X(H)$ satisfy $\e^{(\log\log H)^{\frac {51}{100}}}\le X\le (\log H)^\delta$, and let $1\le m_1\neq m_2\le X$. Set $q\coloneqq \ord{m_2-m_1}$. Let $f_1,f_2\colon \Z\rightarrow \C$ be functions supported on positive integers such that, for some $u>0$, we have \begin{equation}\label{cond_bound}\ord{f_i(n)}\ll u^{\omega(n)}\end{equation} for $i=1,2$.
    Suppose that there exist $\sigma,\tau>0$, $0<\epsilon <1/3$ and a positive measurable function $W_H\colon\R^+\rightarrow\R^+$ such that
    \begin{align}
    \label{cond_unif}
      & \Bigg| \sum_{\substack{n\le T\\n\equiv a \bmod  q}}f_i(n)\Bigg|\ll \frac{q^{\sigma}T^\tau}{W_H(T)}\qquad \text{for every }0\le a<q,\, i\in \{1,2\},\, T>H^{\epsilon}\text{ and}\\
    \label{cond_W}
    & H^{3\epsilon}\ll\int _{H^\epsilon}^T \frac {t^\tau\d t}{W_H(t)}\ll \frac {T^{\tau+1}}{W_H(T)}\qquad \text{for every }T>H^{3\epsilon}    \,.
    \end{align}
     Then, for every $1\le k\le (\log H)^\nu$ and for every $A>0$ we have \begin{equation*}
        \Bigg|\mbf J_{f_1,f_2,k}(m_1,m_2)\Bigg|\ll\frac{H^{d+1}}{k^4(\log H)^A}+\frac{(\log X)^2X^{2(\sigma+\tau d+d+1)}H^{d-1+2\tau}}{k^{4\tau}W\left((d+1)X^dH/k^2\right)^2}\,,
    \end{equation*} where the implied constant depends only on $d,\delta,u,\sigma,\tau,\epsilon,\nu ,A$ and on the constants involved in (\ref{cond_bound}), (\ref{cond_unif}), (\ref{cond_W}).
\end{prop}
\begin{proof} 
    Call  \begin{equation}\label{defS'}S_{f_i}'(k^2\gamma)\coloneqq S_{f_i,k}(\gamma)=\sum_{\substack{h\le (d+1)X^dH\\k^2\mid h}}f_i(h/k^2)e^{i\gamma h}=\sum_{\substack{h'\le \frac{(d+1)X^dH}{k^2}}}f_i(h')e^{ih'(\gamma k^2)}.\end{equation}
    By (\ref{cond_bound}) and by Selberg--Delange's method, we have \begin{equation}\label{bounds}\ord{S_{f_i}'}\ll \sum_{\substack{h'\le \frac{(d+1)X^dH}{k^2}}}\ord{f_i(h')}\ll \frac{X^dH(\log H)^{u-1}}{k^2}\,.\end{equation}
    Let  $\rho=\rho(X,H)\coloneqq\frac{\log X}{H}$ and call $\Gamma\coloneqq (-\pi,\pi]^2$,
    $$\Gamma_\rho^1=\Gamma \cap \{\|\alpha_1-\alpha_2\|_{2\pi}<\rho\}\cap \{\|\alpha_1m_1-\alpha_2m_2\|_{2\pi}<\rho\}\,,\qquad \text{ and }\qquad \Gamma_\rho^2\coloneqq \Gamma\setminus \Gamma_\rho^1\,,$$ where $\|\alpha\|_{2\pi}\coloneqq \min_{r\in \Z}|\alpha-2\pi r|$.
    In analogy to the classical circle method, we call $\Gamma_\rho^1$ the major arcs and $\Gamma_\rho^2$ the minor arcs.
    Then split $\mbf J_{f_1,f_2,k}(m_1,m_2)$ as \begin{align*}&\frac 1 {(2\pi)^2 }\int_{\Gamma_\rho^1}{S'_{f_1}(\alpha_1k^2)}\overline{S'_{f_2}(\alpha_2k^2)}\left(\prod_{j=0}^dK_H(m_1^j\alpha_1-m_2^j\alpha_2)\right)\d \alpha_1\d\alpha_2
    \\+&\frac 1 {(2\pi)^2 }
    \int_{\Gamma_\rho^2}{S'_{f_1}(\alpha_1k^2)}\overline {S'_{f_2}(\alpha_2k^2)}\left(\prod_{j=0}^dK_H(m_1^j\alpha_1-m_2^j\alpha_2)\right)\d \alpha_1\d\alpha_2\,.\end{align*}
    By bounding $\ord {S_{f_i}}$ with (\ref{bounds}) and by bounding $d-1$ heat kernels with Lemma \ref{lemma-HK}(1), the contribution of the minor arcs to $\mbf J(m_1,m_2)$ is bounded by \begin{align*}\ll \frac{\left(X^dH(\log H)^{u-1}\right)^2}{k^4}H^{d-1}\int_{\Gamma_\rho^2}K_H(\alpha_1-\alpha_2)K_H(m_1\alpha_1-m_2\alpha_2)\d\alpha_1\d\alpha_2\,.\end{align*}
    In the inner integral above, change variables as in \cite[Lemma 2.9]{FS} and apply Lemma \ref{lemma-HK}(4--5); then the integral is bounded by $$\int_{\max\{|\beta_1|,|\beta_2|\}\ge \rho}\!\!\!\!\!K_H(\beta_1)K_H(\beta_2)\d\beta_1\d\beta_2\ll \int _{(-\pi,\pi]}\!\!\!K_H(\beta_1)\d\beta_1\int_{|\beta_2|\ge \rho}\!\!\!K_H(\beta_2)\d\beta_2\ll\frac 1 {\rho H\exp\left((\rho H)^2/(4\pi)\right)}\,.$$
    The contribution of the minor arcs to $\mbf J(m_1,m_2)$ is then
    $$\ll\frac {X^{2d}H^{d+1}(\log H)^{2u-2}} {\rho H\exp\left((\rho H)^2/(4\pi)\right)k^4}=\frac {X^{2d}H^{d+1}(\log H)^{2u-2}} {(\log X)\exp\left((\log X)^2/(4\pi)\right)k^4}\ll\frac{H^{d+1}}{(\log H)^Ak^4}$$ for every $A>0$.
    Now we deal with the major arcs. Repeating the same majorizations we used with the minor arcs, the contribution of the major arcs to $\mbf J(m_1,m_2)$ can then be bounded by \begin{align*}&\ll \left\|{S_{f_1,k}(\alpha_1)}\overline {S_{f_2,k}(\alpha_2)}\right\|_\infty^{\Gamma_{\rho}^1}H^{d-1}\!\!\int_{\Gamma_\rho^1}\!K_H(\alpha_1-\alpha_2)K_H(m_1\alpha_1-m_2\alpha_2)\d\alpha_1\d\alpha_2\,,
    \end{align*}
    where the integral is $O(1)$ by Lemma \ref{lemma-HK}(4) and  \cite[Lemma 2.9]{FS}.
     For every $(\alpha_1,\alpha_2)\in \Gamma_\rho^1$ there are two integers $a_1,a_2$ such that \begin{equation}\label{approxalpha}\ord {\alpha_1-\frac{2\pi a_1} q}\ll \frac {\rho m_2} { q}\qquad \text{and}\qquad \ord {\alpha_2-\frac{2\pi a_2} q}\ll \frac {\rho m_1} { q}\,.\end{equation}
    Let $A_i(T)\coloneqq \sum_{t\le T}f_i(t)e^{2\pi i k^2a_it/q}$.
    Assumption (\ref{cond_unif}) implies \begin{equation}\label{boundA1}\ord{A_i(T)}\le \sum_{b\bmod q}\Bigg|\sum_{\substack{t\le T\\t \equiv b\bmod q}}f_i(t)\Bigg|\ll\frac{q^{\sigma+1}T^\tau}{W(T)}\,.\end{equation} Writing $f_i(t)e^{i\alpha_ik^2t} $ as $f_i(t)e^{2\pi ia_ik^2t/q}\cdot e^{ik^2t(\alpha_i-2\pi a_i/q)}$, by partial summation
    \begin{align*}
        \ord{S_{f_i}'(k^2\alpha_i)}&\le \ord{A_i((d+1)X^dH/k^2)}+k^2\ord {\alpha_i-\frac{2\pi a_i} q}\int_1^{(d+1)X^dH/k^2}\ord{A_i(t)}\d t\,.
        \end{align*}
        Using (\ref{boundA1}) and (\ref{approxalpha}), we can bound the second summand on the right-hand side as
        $$\ll \frac{k^2\rho m_{3-i}} q\int_1^{(d+1)X^dH/k^2}\frac{q^{\sigma+1}t^\tau}{W(t)}\d t \ll \frac{\rho  q^\sigma H^{\tau+1}X^{\tau d+d+1}}{k^{2\tau}W((d+1)X^dH/k^2)}\,,$$ where the bound for the integral over $[H^\epsilon,T]$ comes from our hypothesis (\ref{cond_W}), and the same integral over $[1,H^\epsilon]$ is negligible by (\ref{cond_bound}) and (\ref{cond_W}). Thus, using (\ref{boundA1}) to bound the first summand, we obtain
        \begin{align*}
            \ord{S_{f_i,k}'(k^2\alpha_1)}&\ll \frac{{q^{\sigma+1}X^{\tau d}H^\tau}}{k^{2\tau}W\left((d+1)X^dH/k^2\right)}+\frac{{\rho q^\sigma X^{\tau d +d+1}H^{\tau+1}}}{k^{2\tau}W\left((d+1)X^dH/k^2\right)} \ll \frac{(\log X)q^{\sigma}X^{\tau d+d+1}H^\tau}{k^{2\tau}W\left((d+1)X^dH/k^2\right)}\,.
    \end{align*}
    Then, uniformly over $\Gamma_\rho^1$ we have \begin{equation*}
        \left\|{S_{f_1,k}(\alpha_1)}\overline{S_{f_2,k}(\alpha_2)}\right\|^{\Gamma_\rho^1}_\infty\ll \left(\frac{(\log X)X^{\sigma+\tau d+d+1}H^\tau}{k^{2\tau}W\left((d+1)X^dH/k^2\right)}\right)^2\,.\qedhere\end{equation*}
    \end{proof}



\section{Positive multiplicative functions over random polynomials}\label{sec-g3}
For type 1 and R terms, the saving in arithmetic progressions is modest (see \cite[propositions 3.11, 3.17]{Bo2}), so a straightforward application of the machinery of \S\ref{section-secon_moment} would produce an unsatisfactory error term. Thus, in this section we bound the contribution of type 1 and R terms to Theorem \ref{teo_main} by explicitly expanding the second moment. 
We give a detailed proof for type~1 error terms; since the argument is essentially the same, we only sketch the proof for type~R error terms. For type~1 it suffices to recall that $G^{(1)}(D)\ge 0$, together with the following bound from \cite[Equation (3.3)]{Bo2}:
\begin{equation}\label{boundG3}G^{(1)}(D)\ll\mu^2(D)\summ{D'\le C^{13}\\D'\mid D}\leg {13}4^{\omega(D')}\leg 3 4 ^{\omega(D/D')}\,.\end{equation} 
We apply (\ref{boundG3})and estimate the resulting sum using Nair's bound \cite{Na}.

\begin{prop}\label{prop_sm3}
    Fix $d\ge 1$ and $B,\delta,\nu,\eta>0$. Let $H>3$ and let $X$ satisfy $2\le X\le (\log H)^{{\delta}}$. Let $C(H)$ and $z(H)$ be as in (\ref{def_zC}). Then, for every $b\in \{1,2,3,5,6,7\}$ and $1\le k\le (\log H)^\nu$ we have \begin{equation}\label{bound3}
        \sum_{P\in \pol d H}\Bigg|\summ{m\le X\\k^2\mid P(m)\\P(m)/k^2\equiv b\bmod 8}G^{(1)}\bigg(\frac{P(m)}{k^2}\bigg)\Bigg|^2\ll\frac{ X(\log H)^{15}H^{d+1}}{k^2} +\frac{X^2H^{d+1}}{k^2(\log H)^{7/16-377\eta/16}}\,,
    \end{equation} where the implied constant depends only on $d,\delta, \eta,\nu $.
\end{prop}
\begin{proof}
We prove the result for $b\in \{1,3,5,7\}$, but the result can be similarly proven with minor modifications for $b\in \{2,6\}$. First, since $G^{(1)}$ is non-negative, we can drop the condition $P(m)/k^2\equiv b\bmod 8$; this can only increase the second moment.
Split the second moment of $G^{(1)}$ into diagonal term and non-diagonal term. Since $G^{(1)}(D)\ll4^{\omega(D)}$, Proposition \ref{prop-diag1} applies and $$\sum_{P\in \pol d H}\Bigg|\sum_{\substack{m\le X\\k^2\mid P(m)}}G^{(1)}\left(\frac{P(m)}{k^2}\right)^2\Bigg|\ll \frac{X(\log H)^{15}H^{d+1}}{k^2}\,.$$
Concerning the non-diagonal term, we can use (\ref{boundG3}) to write
    \begin{align}
        &\sum_{P\in \pol d H}\Bigg|\sum_{\substack{m_1\neq m_2\le X\\k^2\mid P(m_1),P(m_2)}}G^{(1)}\left(\frac{P(m_1)}{k^2}\right)G^{(1)}\left(\frac{P(m_2)}{k^2}\right)\Bigg|
        \\ &\ll
        \sum_{P\in \pol d H}\!\!\!\!\!\sum_{\substack{m_1\neq m_2\le X\\k^2\mid P(m_1),P(m_2)\\P(m_1)/k^2, P(m_2)/k^2\text{ odd}\sfree}}\sum_{\substack{\mcal D_1'\in \D_{13}(D_1'),\, \mcal D_2'\in \D_{13}(D_2')\\D_1',D_2'\le C^{13}\\D_1'\mid P(m_1)/k^2,\,D_2'\mid P(m_2)/k^2}}
        \!\!\!\!4^{-\omega(D_1')-\omega(D_2')}\leg 3 4^{\omega\leg {P(m_1)}{D_1'k^2}+\omega\leg{P(m_2)}{D_2'k^2}}\notag\\
        &=
        \sum_{\substack{D_1',D_2'\le C^{13}\\D_1',D_2'\text{ odd},\sfree}}\frac{\tau_{13}(D_1')\tau_{13}(D_2')}{4^{\omega(D_1')+\omega(D_2')}}\sum_{P\in \pol d H}\sum_{\substack{m_1\neq m_2\le X\\k^2D_i'\mid P(m_i)\\P(m_i)/k^2 \text{ odd,}\sfree}}\leg 3 4^{\omega\leg {P(m_1)}{D_1'k^2}+\omega\leg{P(m_2)}{D_2'k^2}}\,.\label{eq2}
    \end{align}
    Writing $P(t)=c_dt^d+\ldots+c_1t+c_0$ and setting $N_i\coloneqq c_dm_i^d+\ldots+c_1m_i$ for $i=1,2$, we have $P(m_i)=N_i+c_0$, and the condition $k^2D_i'\mid P(m_i)$ can be written as $c_0\equiv -N_i\bmod k^2D_i'$. The sum inside (\ref{eq2}) is then \begin{equation}\label{eq1}\ll\sum_{c_d,\ldots,c_1\in [0,H]}\sum_{m_1\neq m_2\le X}\sum_{\substack{c_0\in [0,H]\\(c_0+N_i)/k^2 \sfree \\c_0\equiv -N_i\bmod k^2D_i'}}\leg 3 4^{\omega\leg {c_0+N_1}{k^2D_1'}+\omega\leg{c_0+N_2}{k^2D_2'}}\,.\end{equation}

    Letting $d_0=\gcd(D_1',D_2')$, the two congruences $c_0\equiv -N_i\bmod k^2D_i'$ have a solution if and only if $k^2d_0\mid (N_1-N_2)$; if so, let $N_0$ be the smallest positive integer satisfying $N_0\equiv N_i\bmod k^2D_i'$ for $i=1,2$.
    Then let $e_i=\gcd\left(\frac{N_2-N_1}{k^2},D_i'\right)/d_0$ and let $d_i=D_i'/d_0e_i$.
    In this way, $D_i'=d_0d_ie_i$ and, since $D'_1,D'_2$ are square-free, $d_0,d_1,d_2,e_1,e_2$ are also square-free and pairwise coprime. Write any $c_0\in [0,H]$ satisfying $c_0\equiv -N_i\bmod k^2D_i'$ as $c_0=-N_0+c_0'\cdot d_0d_1d_2e_1e_2k^2$ for $c_0'\in [0,H/d_0d_1d_2e_1e_2k^2]$.
    With this choice, \begin{align*}\frac{P(m_1)}{k^2D_1'}&=\frac{N_1-N_0}{k^2D_1'}+c_0'd_2e_2=e_2\left(\frac{N_1-N_0}{k^2D_1'e_2}+c_0'd_2\right)\,\text{ and }\\\frac{P(m_2)}{k^2D_2'}&=\frac{N_2-N_0}{k^2D_2'}+c_0'd_1e_1=e_1\left(\frac{N_2-N_0}{k^2D_2'e_1}+c_0'd_1\right)\,,\end{align*}
    where $\frac{N_1-N_0}{D_1'e_2k^2}$ and $\frac{N_2-N_0}{D_2'e_1k^2}$ are both positive integers.
Define $Q_1,Q_2\in \Z[t]$ by $$Q_1(t)\coloneqq \frac{N_1-N_0}{k^2D_1'e_2}+d_2t\qquad\text{and}\qquad Q_2(t)\coloneqq \frac{N_2-N_0}{k^2D_2'e_1}+d_1t\,,$$ so that $P(m_1)=k^2D'_1e_2Q_1(c_0')$ and $P(m_2)=k^2D'_2e_1Q_2(c_0')$. We replace the condition \lq\lq$P(m_i)/k^2$ square-free\rq\rq\ in (\ref{eq1}) with the weaker condition \lq\lq$Q_i(c_0')\sfree$\rq\rq; in this way we can bound (\ref{eq1}) by
\begin{equation}
    \label{eq4}\ll\sum_{\substack{c_1,\ldots,c_d\in [0,H]\\m_1\neq m_2\le X\\d_0e_1e_2k^2\mid N_2-N_1\\\gcd\left(\frac{N_2-N_1}{k^2},d_1d_2\right)=1}}\leg 3 4^{\omega(e_1)+\omega(e_2)}\sum_{\substack{c_0'\in [0,H/d_0d_1d_2e_1e_2k^2]\\Q_i(c_0') \sfree}}\leg 3 4^{\omega(Q_1(c_0'))+\omega(Q_2(c_0'))}\,.
\end{equation}
Applying the Cauchy--Schwarz inequality to the second sum leads to \begin{equation}\label{cs-appl}\ll \left(\sum_{c_0'\in [0,H/d_0d_1d_2e_1e_2k^2]}\!\!\leg 9{16}^{\omega(Q_1(c_0'))}\!\!\mu^2(Q_1(c_0'))\!\right)^{1/2}\left(\sum_{c_0'\in [0,H/d_0d_1d_2e_1e_2k^2]}\!\!\leg 9{16}^{\omega(Q_2(c_0'))}\!\!\mu^2(Q_2(c_0'))\!\right)^{1/2}\,.\end{equation}
By construction, the coefficients of each $Q_i$ are pairwise coprime, so by Lemma \ref{lemma_nair} applied to the two sums above, we can bound (\ref{eq4}) by 
$$\ll \!\frac{ H}{k^2d_0d_1d_2e_1e_2(\log H)^{7/16}} \!\!\!\!\!\!\!\!\!\sum_{\substack{c_1,\ldots,c_d\in [0,H]\\m_1\neq m_2\le X\\d_0e_1e_2k^2\mid N_2-N_1\\\gcd\left(\frac{N_2-N_1}{k^2},d_1d_2\right)=1}}\!\!\!\!\!\!\!\!\leg 3 4^{\omega(e_1)+\omega(e_2)}\!\!2^{\omega(d_1)+\omega(d_2)}
\!\ll\! \frac{X^2H^{d+1}}{k^2(\log H)^{7/16}}\frac{d_02^{\omega(D_1'/d_0)+\omega(D_2'/d_0)}}{D_1'D_2'}.$$ Here we have bounded $\frac 3 4$ by $1$ and the constants $C_{9/16,d_i}$ by $2^{\omega(d_i)}\le 2^{\omega(D_i'/d_0)}$.
 The quantity  (\ref{eq2}), is then bounded by \begin{align*}
    \ll\frac{X^2H^{d+1}}{k^2(\log H)^{7/16}}\sum_{\substack{D_1',D_2'\le C^{13}}}\frac{\tau_{13}(D_1')\tau_{13}(D_2')\mu^2(D_1')\mu^2(D_2')d_0 2^{\omega(D_1'/d_0)}2^{\omega(D_2'/d_0)}}{D_1'D_2'4^{\omega(D_1')+\omega(D_2')}}
\end{align*}
Writing $D'_i=d_0D''_i$, and dropping the coprimality condition between $d_0,D_i''$ and the modulo 2 conditions, the quantity in (\ref{eq2}) can be bounded by \begin{align*}
    &\ll\sum_{d_0\le C^{13}}\frac{\mu^2(d_0)\tau_{13}(d_0)^2}{d_016^{\omega(d_0)}}\sum_{D_1'',D_2''\le C^{13}/d_0}\frac{\tau_{26}(D_1'')\mu^2(D_1'')}{D_1''4^{\omega(D_1'')}}\frac{\tau_{26}(D_2'')\mu^2(D_2'')}{D_2''4^{\omega(D_2'')}}
    \\&
    \ll\sum_{d_0\le C^{13}}\frac{\mu^2(d_0)\tau_{13}(d_0)^2}{d_016^{\omega(d_0)}}\left(\log\left(\frac{2C^{13}}{d_0}\right)\right)^{13}\ll(\log C(H))^{377/16}\ll(\log H)^{377\eta/16}\,.
\end{align*}
\end{proof}
\begin{remark}
    We used the Cauchy--Schwarz inequality in (\ref{cs-appl}) to simplify the proof, but more refined versions of Nair's bound \cite{He} (or, alternatively, Hooley neutralizers) lead to a factor $(\log H)^{-1/2+\epsilon}$ instead of $(\log H)^{-7/16+\epsilon}$ in (\ref{bound3}).
\end{remark}

 For type R terms  we know from \cite{Bo2} that $G^{(\mathrm{R}) }(D)\ge 0$ and 
 \begin{equation}\label{bound_GR}G^{(\text{R})}(D)\ll \mu^2(D)\summ{D'\le C\\D'\mid D}\leg 3 4 ^{\omega\left(D/D'\right)}\,.\end{equation}
\begin{prop}\label{prop_smR}
    Fix $d\ge 1$ and $\delta,\nu,\eta>0$. Let $H>3$ and let $X$ satisfy $2\le X\le (\log H)^{{\delta}}$. Let $C(H)$ be as in (\ref{def_zC}). Then, for every $b\in \{1,2,3,5,6,7\}$ and $1\le k\le (\log H)^\nu$ we have \begin{equation}\label{boundR}
        \sum_{P\in \pol d H}\Bigg|\summ{m\le X\\k^2\mid P(m)\\P(m)/k^2\equiv b\bmod 8}G^{(\mathrm{R})}\bigg(\frac{P(m)}{k^2}\bigg)\Bigg|^2\ll\frac{ X(\log H)^{15}H^{d+1}}{k^2} +\frac{X^2H^{d+1}}{k^2(\log H)^{7/16-377\eta/16}}\,,
    \end{equation} where the implied constant depends only on $d,\delta, \eta,\nu $.
\end{prop}
\begin{proof}
    Just replace (\ref{boundG3}) with (\ref{bound_GR}) in the proof of Proposition \ref{prop_sm3}.  
   
\end{proof}


\section{Sieving square-full parts}\label{sec-ssf}
In this section, we use a version of the square-free sieve to show that Theorem \ref{teo_mainES} follows from Theorem \ref{teo_main}. Recall that the $2$-Selmer rank $s(n)$ is independent on the square-full part of $n$.
We isolate the largest square divisor $k^2$ of $P(m)$:  when $k$ is large Nair's bound suffices, whereas for small $k$ we can use Theorem \ref{teo_main}.
\begin{proof}[Proof of Theorem \ref{teo_mainES}]
    First, we use the relation $s(D)=s(k^2D)$ to split the desired sum as \begin{align*}
    \ll &
    \sum_{b'\in \{1,2,3,5,6,7\}}\sum_{P\in \pol d H}\Bigg|\sum_{\substack{1\le k\ll \sqrt{HX^d}\\k^2b'\equiv b\bmod 8}}\sum_{\substack{m\le X\\k^2\mid P(m)\\ P(m)/k^2\equiv b'\bmod 8\\P(m)/k^2\sfree}}\!\!\!\!\!\!\!\!\left(2^{s\left(\frac{P(m)}{k^2}\right)}-3\right)\Bigg|^2\,.
    \end{align*} 
    The advantage is that $\frac{P(m)}{k^2}$ is square-free and belongs to a fixed residue class modulo $8$. Fix $b'\in \{1,2,3,5,6,7\}$.
    Letting $\gamma(H)=(\log H)^{1000}$, we split the sum above as
    \begin{align*}
    \ll
    \sum_{P\in \pol d H}\Bigg|\!\sum_{\substack{1\le k\le\gamma\\k^2b'\equiv b\bmod 8}}\!\sum_{\substack{m\le X\\k^2\mid P(m)\\P(m)/k^2\equiv b'\bmod 8\\P(m)/k^2\sfree}}\!\!\!\!\!\!\!\!\!\!\!\!\left(2^{s\left(\frac{P(m)}{k^2}\right)}-3\right)\Bigg|^2
    \!+\!\!
    \sum_{P\in \pol d H}\!\!\Bigg(\sum_{\substack{\gamma< k\ll \sqrt{H X^d}}}\,\!\!\!\!\!\!\!\!\!\!\sum_{\substack{m\le X\\k^2\mid P(m)\\P(m)/k^2\sfree}}\!\!\!\!\!\!\!\!\!\!\!\left|2^{s\left(\frac{P(m)}{k^2}\right)}-3\right|\!\Bigg)^2\!\!.\end{align*}
First, we bound the sum over large values of $k$. Expand the second moment in the sum above. From (\ref{bound_sel}) it follows $\ord{2^{s(D)}-3}\le 9\cdot 4^{\omega(D)}$, so we have \begin{align*}
    \ll\sum_{m_1,m_2\le X}\sum_{\substack{\gamma<k_1,k_2\ll\sqrt{H X^d}}}\sum_{\substack{P\in\pol d H\\k_i^2|P(m_i)\\P(m_i)/k_i^2\sfree}}{4^{\omega\left(\frac{P(m_1)}{k_1^2}\right)}}{4^{\omega\left(\frac{P(m_2)}{k_2^2}\right)}}\,.
\end{align*}
In the sum above, use $\omega(P(m_2)/k_2^2)\le\omega(P(m_2))$ and note that there is at most one $\gamma<k_2\ll\sqrt{H X^d}$ with $k_2^2\mid P(m_2)$ and $P(m_2)/k_2^2$ square-free; thus, we have
 \begin{align*}
    &\ll\sum_{m_1,m_2\le X}\sum_{\substack{\gamma<k_1\ll\sqrt{H X^d}}}\sum_{\substack{P\in\pol d H\\k_1^2|P(m_1)}}{4^{\omega\left(\frac{P(m_1)}{k_1^2}\right)}}{4^{\omega\left({P(m_2)}\right)}}
    \end{align*}
    Fix the values of $m_1$ and $m_2$, and split the sum over $k_1$ as 
    $$\sum_{\substack{\gamma<k_1\le {H^{0.2}}}}\sum_{\substack{P\in\pol d H\\k_1^2|P(m_1)}}{4^{\omega\left(\frac{P(m_1)}{k_1^2}\right)}}{4^{\omega\left({P(m_2)}\right)}}\,
    + \sum_{\substack{H^{0.2}<k_1\ll\sqrt{H X^d}}}\sum_{\substack{P\in\pol d H\\k_1^2|P(m_1)}}{4^{\omega\left(\frac{P(m_1)}{k_1^2}\right)}}{4^{\omega\left({P(m_2)}\right)}}\,.$$
    The sum in the larger range can be bounded by 
    \begin{align*}
        \ll\sum_{\substack{H^{0.2}<k_1\ll\sqrt{H X^d}}}\sum_{\substack{P\in\pol d H\\k_1^2|P(m_1)}}H^{0.1}= \sum_{\substack{H^{0.2}<k_1\ll\sqrt{H X^d}}}\left(\frac{H^{d+1.1}}{k_1^2}+O\left(H^{d+0.1}\right)\right)\ll H^{d+0.9}\,,
    \end{align*} and, after summing over $m_1,m_2$, this  gives a contribution of $\ll X^2H^{d+0.9}$.
Concerning the sum over $k_1$ in the intermediate range, we can apply H\"older's inequality with exponents $3/2$ and $3$ to bound it by
    \begin{align}\label{eq6}
    \sum_{\gamma<k_1\le H^{0.2}}\Bigg(\sum_{\substack{P\in\pol d H\\k_1^2|P(m_1)}}8^{\omega\left(\frac{P(m_1)}{k_1^2}\right)}\Bigg)^{2/3}\Bigg(\sum_{P\in \pol d H}64^{\omega(P(m_2))}\Bigg)^{1/3}\,.
\end{align}
Let $P(m_1)=c_dm_1^d+\ldots+c_1m_1+c_0\coloneqq N_1+c_0$. After rewriting the sum inside the parentheses on the left, we can bound it with Lemma \ref{lemma_nair} as
$$\sum_{c_1,\ldots, c_d\in (0,H]} \sum_{n\in \left[\frac{N_1}{k_1^2},\frac{N_1+H}{k_1^2}\right]\cap \Z}8^{\omega\left(n\right)}\ll \frac{H^{d+1}(\log H)^7}{k_1^2}\,. $$ Observe that the condition $k_1\le H^{0.2}$ is crucial for the applicability of Lemma \ref{lemma_nair}. In a similar way, we can bound the sum inside the parentheses on the right in (\ref{eq6}) by $\ll H^{d+1}(\log H)^{63}$, so (\ref{eq6}) is bounded by
$$\ll H^{d+1}(\log H)^{77/3}\sum_{\gamma<k_1\le H^{0.2}}k_1^{-4/3}\ll \frac{H^{d+1}(\log H)^{77/3}}{(\log H)^{1000/3}}\,.$$

We now move to the sum over small values of $k$. Recalling the definition (\ref{defpsi}), we have \begin{align*}
    \sum_{P\in \pol d H}\bigg|\sum_{\substack{1\le k\le \gamma\\k^2b'\equiv b\bmod 8}}&\Psi_{P,k,b'}(X)\bigg|^2
    \le
    \sum_{P\in \pol d H}\sum_{\substack{1\le k_1,k_2\le \gamma\\k_1^2b',k_2^2b'\equiv b\bmod 8}}\bigg|\Psi_{P,k_1,b'}(X)\Psi_{P,k_2,b'}(X)\bigg|
    \\&\le
    \sum_{\substack{1\le k_1,k_2\le \gamma\\k_1^2b',k_2^2b'\equiv b\bmod 8}}\!\!\!\Bigg(\sum_{P\in \pol d H}\!\!\Psi_{P,k_1,b'}(X)^2\Bigg)^{1/2}\Bigg(\sum_{P\in \pol d H}\Psi_{P,k_2,b'}(X)^2\Bigg)^{1/2}
    \\&\le
    \left(\sum_{\substack{1\le k\le \gamma\\k^2b'\equiv b\bmod 8}}\Bigg(\sum_{P\in \pol d H}\Psi_{P,k,b'}(X)^2\Bigg)^{1/2}\right)^2\end{align*}
By Theorem \ref{teo_main} the quantity above can be bounded by \begin{equation*}
\ll \frac{X^2H^{d+1}}{(\log H)^{7/16-\epsilon/2}}\left(\sum_{\substack{1\le  k\le (\log H)^{1000}\\k^2b'\equiv b\bmod 8}}\frac 1 k\right)^2\ll \frac{X^2H^{d+1}(\log\log H)^{2}}{(\log H)^{7/16-\epsilon/2}} \ll \frac{X^2H^{d+1}}{(\log H)^{7/16-\epsilon}} \,.\qedhere\end{equation*}
\end{proof}



\bibliographystyle{plainurl}
\bibliography{sample}

\end{document}